\documentclass[12pt]{amsart}
\usepackage{amsfonts,amssymb,amsmath,amscd,amsthm}
\usepackage{stmaryrd,graphicx}
\usepackage[all]{xy}

\newtheorem{thm}{Theorem}
\newtheorem{lem}[thm]{Lemma}
\newtheorem{prop}[thm]{Proposition}

\newtheorem{cor}[thm]{Corollary}

\theoremstyle{definition}
\newtheorem{dfn}[thm]{Definition} 
\newtheorem{ex}[thm]{Example}
\newtheorem{rmk}[thm]{Remark}

\numberwithin{thm}{section}
\numberwithin{equation}{section}

\newcommand{\Proof}{\noindent {\it Proof}.\ \ }

\newcommand{\Ext}{\operatorname{Ext}}

\newcommand{\Hilb}{\operatorname{Hilb}}
\newcommand{\HF}{\operatorname{HF}}
\newcommand{\red}{\operatorname{red}}

\newcommand{\Bs}{\operatorname{Bs}}

\newcommand{\Pic}{\operatorname{Pic}}

\newcommand{\coker}{\operatorname{coker}}

\newcommand{\car}{\operatorname{char}}

\newcommand{\codim}[1]{\operatorname{codim}}

\renewcommand{\labelenumi}{{\rm (\arabic{enumi})}}

\title{Obstructions to deforming curves on a prime Fano $3$-fold}
\author{Hirokazu Nasu}
\date{}
\address{
Department of Mathematical Sciences,
Tokai University,
4-1-1 Kitakaname, Hiratsuka, 
Kanagawa 259-1292, JAPAN}
\email{nasu@tokai-u.jp}

\begin{document}


\begin{abstract}
We prove that for every smooth prime Fano $3$-fold $V$,
the Hilbert scheme $\Hilb^{sc} V$ of smooth connected curves on $V$
contains a generically non-reduced irreducible component
of Mumford type.
We also study the deformations of degenerate curves $C$ in $V$, 
i.e., curves $C$ contained in a smooth anticanonical member 
$S \in |-K_V|$ of $V$.
We give a sufficient condition for $C$ to be stably degenerate, 
i.e., every small (and global) deformation of $C$ in $V$ is contained 
in a deformation of $S$ in $V$.
As a result, by using the Hilbert-flag scheme of $V$,
we determine the dimension and the smoothness
of $\Hilb^{sc} V$ at the point $[C]$, assuming that
the class of $C$ in $\Pic S$ is generated 
by $\mathbf h:=-K_V\big{\vert}_S$ 
together with the class of a line, or a conic on $V$.
\end{abstract}
\maketitle

\section{Introduction}
\label{sect:introduction}

We work over an algebraically closed field $k$ of characteristic $0$.
Let $X$ be a smooth Fano $3$-fold of index $r$ 
with Picard group $\Pic X\simeq \mathbb Z$.
Then by \cite{Iskovskih77,Iskovskih78},
all such $X$ are classified into $17$ classes up to deformation equivalence
and we have $1 \le r \le 4$.
Let $\Hilb^{sc} X$ denote the Hilbert scheme  of smooth connected curves in $X$.
Mumford~\cite{Mumford} first 
proved that if $r=4$ (i.e.~$X\simeq \mathbb P^3$), then
$\Hilb^{sc} X$ contains a generically non-reduced (irreducible) component.
This example was generalized in \cite{Gruson-Peskine,Kleppe87,Ellia87}, etc.,
for $X=\mathbb P^3$, and in \cite{Mukai-Nasu,Nasu4} 
for many uniruled $3$-folds $X$.
It is known that if $r=3$ (i.e.~$X \simeq Q^3 \subset \mathbb P^4$)
or $r=2$ (i.e.~$X$ is a {\em del Pezzo} $3$-fold), 
then $\Hilb^{sc} X$ contains (infinitely many)
generically non-reduced components.

In this paper, we discuss the existence of
a generically non-reduced component of $\Hilb^{sc} X$
for every $X$ with $r=1$, i.e., a {\em prime} Fano $3$-fold $X$,
in the view point of a further generalization of Mumford's example.
Let $V$ be a prime Fano $3$-fold of genus $g$ $(:=(-K_V)^3/2+1)$,
and let $\Hilb^{sc}_{d,p} V$ denote
the subscheme of $\Hilb^{sc} V$ parametrising curves of degree $d$
and genus $p$. 
\begin{thm}
  \label{thm:main1}
  The Hilbert scheme $\Hilb^{sc}_{4g,4g+1} V$
  contains a generically non-reduced irreducible
  component (of Mumford type) of dimension $5g+1$,
  whose general member $C$ satisfies:
  \begin{enumerate}
    \item $C$ is contained in a smooth anticanonical member $S \in |-K_V|$,
    \label{item:degeneracy}
    \item $C$ belongs to the class $-2K_V\big{\vert}_S+2E$ in $\Pic S$
    for a good conic $E\simeq \mathbb P^1$ on $V$,
    i.e., a conic with trivial normal bundle 
    $N_{E/V} \simeq \mathcal O_{\mathbb P^1}^2$, and
    \label{item:linear equivalence}
    \item $h^0(C,N_{C/V})=5g+2$.
    \label{item:tangential dimension}
  \end{enumerate}
\end{thm}
This is a generalization of a result in \cite{Nasu5}
(for $g=3$).
We will sketch its proof.
As is well known, $V$ contains good conics $E$, which
are parametrised by an open subset $\Gamma'$ of
the Fano surface $\Gamma:=\Hilb_{2,0} V$ of $V$
(cf.~Lemma~\ref{lem:good conic}).
The pairs $(E,S)$ of $E$ and 
a smooth member $S \in |-K_V|$ containing $E$
are parametrised by an open subset $U$ of a 
$\mathbb P^{g-2}$-bundle over $\Gamma'$. (Thus $\dim U=g$.)
We consider the maximal family $W$ of curves $C$ 
contained in a smooth $S \in |-K_V|$, and 
belonging to the linear equivalence class in \eqref{item:linear equivalence}
for some $E$. Then $W$ is isomorphic to an open subset of 
a $\mathbb P^{4g+1}$-bundle over $U$,
and thus $\dim W=5g+1$.
(See \S\ref{subsec:construction} for the construction of $W$.)
On the other hand, 
we compute that $h^0(C,N_{C/V})=5g+2$, using 
the fact that the Hilbert-flag scheme $\HF^{sc} V$
is nonsingular 
at $(C,S)$ of expected dimension $5g+1$ 
(cf.~Lemma~\ref{lem:non-lifting of curves}).
By using a result in \cite{Nasu5} (cf.~Lemma~\ref{lem:k3 and fano}), 
we show that every $C$ is obstructed in $V$.
Then as in Mumford's example,
there exists an inequality 
$\dim W \le \dim_{[C]} \Hilb^{sc} V < h^0(N_{C/V})=\dim W+1$,
and this inequality immediately implies that
the closure $\overline W$ of $W$ in $\Hilb^{sc} V$
is an irreducible component of $(\Hilb^{sc} V)_{\red}$, 
and $\Hilb^{sc} V$ is generically non-reduced along $\overline W$.
In Table~\ref{table:GNRC of Mumford type},
for every integer $1 \le r \le 4$
and every smooth Fano $3$-fold $X$ of index $r$ with 
$\Pic X \simeq \mathbb Z$,
we list a series of generically non-reduced components 
$W \subset \Hilb^{sc} X$ ``of Mumford type'':
\begin{table}[h]
  \caption{Generically non-reduced components of Mumford type}
  \label{table:GNRC of Mumford type}
  \begin{minipage}{13cm}
    \begin{center}
      \begin{tabular}{|c|c|c|c|c|l|}
	\hline
	$r$ & class of $S$ & class of $C$ & $E$ & $\dim W$ & \\
	\hline
	$4$ & $-\frac 34 K_X$ & & & $56$ & Mumford~\cite{Mumford}
	(cf.~Ex.~\ref{ex:mumford})\\
	\cline{1-2}\cline{5-6} 
	$3$ & $-\frac 23 K_X$ & $-K_X\big{\vert}_S+2E$ & line & $42$ & \cite{Mukai-Nasu} (cf.~Ex.~\ref{ex:non-reduced index 3})\\
	\cline{1-2}\cline{5-6} 
	$2$ & $-\frac 12 K_X$ & & & $4n+4$\footnote{
	  Here $n$ denotes the degree of a del Pezzo $3$-fold $X$
	  (i.e.~$n=(-K_X)^3/8$).} 
	& \cite{Mukai-Nasu,Nasu4}\\
	\hline
	$1$ & $-K_X$ & $-2K_X\big{\vert}_S+2E$ & conic & $5g+1$ & 
	\cite{Nasu5} ($g=3$)\\
	\hline
      \end{tabular}
    \end{center}
  \end{minipage}
\end{table}

As the second topic of this paper,
we study the deformations of degenerate curves $C$ in $V$,
i.e., curves contained in some member $S\in |-K_V|$.
We are interested in
(i) the stability of the degeneration of $C$,
and also (ii) the (un)obstructedness of $C$ in $V$.
Here we say $C$ is {\em stably degenerate}
if every small and global deformation $C'$ of $C$ in $V$ is contained 
in some deformation $S'$ of $S$ in $V$
(cf.~Definition~\ref{dfn:stably degenerate}).
\begin{thm}
  \label{thm:main2}
  Let $C$ be a smooth connected curve of genus $g(C)$ on $V$
  contained in a smooth anticanonical member $S\in |-K_V|$ of $V$,
  and let $\mathbf h$ denote the class of $-K_V\big{\vert}_S$ in $\Pic S$.
  Suppose that we have either
  \begin{enumerate}
    \renewcommand{\theenumi}{{\roman{enumi}}}
    \renewcommand{\labelenumi}{{\rm [\theenumi]}}
    \item $C \sim n \mathbf h$ for some integer $n$,
    \label{item:c.i.}
    \item $S$ contains a good line $E$ 
    (cf.~Example~\ref{ex:obstruction of flag including lines}) 
    on $V$
    and $C \sim a \mathbf h + b E$ for some integers $a,b\ge 0$
    with $(a,b)\ne (0,1)$, or
    \label{item:w.line}
    \item $S$ contains a good conic $E$ on $V$ and
    $C \sim a \mathbf h + b E$ for some integers $a,b\ge 0$
    with $(a,b)\ne (0,1)$.
    \label{item:w.conic}
  \end{enumerate}
  Then 
  \begin{enumerate}
    \item $C$ is stably degenerate in $V$.
    \item 
    $
    \dim_{[C]} \Hilb^{sc} V=
    \begin{cases}
      g+g(C)+1 & \mbox{if [\ref{item:c.i.}] with $n\ge 2$ holds, and} \\
      g+g(C) & \mbox{otherwise.} \\
    \end{cases}
    $
    \item $C$ is obstructed in $V$ 
    if and only if [\ref{item:w.conic}] with $a=b\ge 2$ holds.
  \end{enumerate}
\end{thm}
The same subject was studied in \cite{Nasu4}
for degenerate curves on a del Pezzo $3$-fold.
By~\cite{Moisezon}, if $S$ is general in $|-K_V|$, 
then every smooth curve $C$ on $S$ is of type [\ref{item:c.i.}].
We call a curve $C$ on $S$ of this type 
a {\em complete intersection in $S$}.
On the other hand, in the studies of Fano $3$-folds $V$
(cf.~\cite{Iskovskih77,Iskovskih78}), 
lines and conics on $V$ play an important role.
Thus, the curves on $S$ of type [\ref{item:w.line}] or [\ref{item:w.conic}]
seem to be of the next natural class
to consider their deformations in $V$,
other than the complete intersections in $S$.
In Propositions~\ref{prop:complete intersection},
\ref{prop:with line} and \ref{prop:with conic}, 
we give the dimension of $\Hilb^{sc} V$ at $[C]$ more explicitly, 
in terms of $(a,b)$ and $n$.

The organization of this paper is as follows.
In \S\ref{subsec:prime fano}, we recall some basic results
on prime Fano $3$-folds.
In \S\ref{subsec:flag scheme}, we discuss the Hilbert-flag schemes.
We consider the image of the first projection
$pr_1: \HF^{sc} V \rightarrow \Hilb^{sc} V$ sending $(C,S)$ to $[C]$,
and prove Theorem~\ref{thm:codimension of flag in hilb},
which is a generalization of 
a result due to Kleppe~\cite{Kleppe87} for $V=\mathbb P^3$,
although its proof is not new.
In \S\ref{subsec:curve and k3}, we apply this theorem
to our case: $C \subset S \subset V$, where 
$C$ is a curve, $S$ is a $K3$ surface,
and $V$ is a Fano $3$-fold
(cf.~Lemma~\ref{lem:codimension of S-maximal in hilb}),
and prepare the two key lemmas
(cf.~Lemmas~\ref{lem:non-lifting of curves} and \ref{lem:k3 and fano}).
We prove Theorems~\ref{thm:main1} in \S\ref{sect:non-reduced}
and \ref{thm:main2} in \S\ref{sect:stably degenerate}, respectively.

\section{Preliminaries}
\label{sect:prelimiaries}

\subsection{Prime Fano $3$-folds}
\label{subsec:prime fano}

In this section, we recall some basic facts on prime Fano $3$-folds.
We refer to Iskovskih~\cite{Iskovskih78}, 
or a survey in \cite{AG5} for the details.
A smooth projective $3$-fold $V$ is called a 
{\em Fano $3$-fold} if the anticanonical divisor $-K_V$ of $V$ is ample.
The maximal integer $r$ such that $-K_V \sim rH$
for some Cartier divisor $H$ is called the 
{\em (Fano) index} of $V$, and we have $1 \le r\le 4$ for every $V$.
We consider a {\em prime} Fano $3$-fold $V$, i.e., a Fano $3$-fold $V$
with $r=1$, and such that $\Pic V$ is generated by $H:=-K_V$.
Then the genus $g:=H^3/2+1$ of $V$
can be any integer between $2$ and $12$, except $11$.
The linear system $|H|$ on $V$ defines a morphism
$\Phi_{|H|}: V \rightarrow \mathbb P^{g+1}$,
which is an embedding, 
or a finite morphism of degree $2$ onto its image in $\mathbb P^{g+1}$.
In the latter case, $V$ is said to be {\em hyperelliptic}.
By virtue of the classification due to Iskovskih, together with
Mukai's work (\cite{Mukai88,Mukai92a}),
in which the anticanonical models of $V$
(i.e.~$\Phi_{|H|}(V) \subset \mathbb P^{g+1}$) were described
as linear or quadratic sections of homogeneous spaces
$\Sigma$ ($g=7,9,10$),
every prime Fano $3$-fold $V$ of genus $g$
is isomorphic to $V_{2g-2}$ ($2 \le g\le 12$) or $V_4'$
in Table~\ref{table:prime fano 3-folds} (cf.~\cite{Mukai95}).
\begin{table}[h]
  \caption{Prime Fano $3$-folds}
  \label{table:prime fano 3-folds}
  \begin{minipage}{16cm}
    \begin{center}
      \begin{tabular}{|c|l|}
	\hline
	$g$ & anticanonical model (or morphism)\\
	\hline
	$2$ & $V_2 \overset{2:1}{\longrightarrow} \mathbb P^3$:
	a double cover branched along $(6) \subset \mathbb P^3$\\
	\hline
	$3$ & $V_4 = (4) \subset \mathbb P^4$: a quartic hypersurface\\
	\hline
	& $V_4' \overset{2:1}{\longrightarrow} (2) \subset \mathbb P^4$:
	a double cover branched along $(2)\cap (4) \subset \mathbb P^4$\\
	\hline
	$4$ & $V_6 = (2) \cap (3) \subset \mathbb P^5$: a complete intersection
	of a quadric and a cubic\\
	\hline
	$5$ & $V_8 = (2) \cap (2) \cap (2) \subset \mathbb P^6$:
	a complete intersection of three quadrics\\
	\hline
	$6$ & $V_{10} = [V_5\footnote{
	    $V_5$ is a del Pezzo $4$-fold 
	    $V_5=[G(2,5) \subset \mathbb P^9] \cap \mathbb P^7$,
	    or a cone over a quintic del Pezzo $3$-fold.}
	  \subset \mathbb P^7]\cap (2)$: 
	a quadratic hypersurface section of a del Pezzo $4$-fold\\
	\hline
	$7$ & $V_{12} = [\Sigma^{10}_{12}=SO(10)/U(5) \subset \mathbb P^{15}]\cap \mathbb P^8$:
	a linear section of a orthogonal Grassmannian \\
	\hline 
	$8$ & $V_{14} = [G(2,6) \subset \mathbb P^{14}]
	\cap \mathbb P^9$: a linear section of a Grassmannian \\
	\hline
	$9$ & $V_{16} = [\Sigma^{6}_{16}=Sp_6(6)/U(3) \subset \mathbb P^{13}] \
	\cap \mathbb P^{10}$: a linear section of a symplectic Grassmannian \\
	\hline
	$10$ & $V_{18} = [\Sigma^{5}_{18} \subset \mathbb P^{13}] \cap \mathbb P^{11}$: a linear section of a $G_2$-variety \\
	\hline
	$12$ & $V_{22}\footnote{
	  $V_{22}$ is isomorphic to the variety
	  $G(3,7,N) \subset \mathbb P^{13}$ associated with a 
	  non-degenerate $3$-dimensional subspace 
	  $N \subset \wedge^2 k^7$.}
	\subset \mathbb P^{13}$: a Mukai-Umemura $3$-fold
	(cf.~\cite{Mukai-Umemura83})\\
	\hline
      \end{tabular}
    \end{center}
  \end{minipage}
\end{table}
Every general member of $|-K_V|$ is a smooth $K3$ surface
(cf.~\cite{Shokurov79b}), and a smooth projective
$3$-fold $X \subset \mathbb P^{g+1}$
(for $g \ge 3$) is a prime Fano $3$-fold (of genus $g$)
if every general linear section
$[X \subset \mathbb P^{g+1}] \cap \mathbb P^{g-1}$ of codimension $2$
is a canonical curve (of genus $g$ in $\mathbb P^{g-1}$).

We next recall the geometry of 
lines and conics on prime Fano $3$-folds.
By a {\em line} (resp. a {\em conic}) on $V$, 
we mean a reduced irreducible rational curve $E$ on $V$
with $(E.H)_V=1$ (resp. $(E.H)_V=2$).
There exists a line (cf.~\cite{Shokurov79a}) and a conic on every $V$, and
$V$ is not covered by the family of lines (because its dimension is
at most $1$), but that of conics.
Let $E$ be a conic on $V$.
Then since $V$ does not contain a plane,
by e.g., \cite[Lemma~4.2 and Proposition~4.3]{Iskovskih78},
$E$ is one of the following type:
\begin{equation}
  \label{eqn:type of conic}
  (k,-k):
  \qquad
  N_{E/V} \simeq \mathcal O_{\mathbb P^1}(k) 
  \oplus \mathcal O_{\mathbb P^1}(-k),
  \qquad
  k=0,1,2.
\end{equation}
In this paper, if $N_{E/V}$ is trivial, $E$ is called a {\em good conic}, 
and a {\em bad conic} otherwise.
For every $V$, the Hilbert scheme $\Gamma$ ($=\Hilb_{2,0}^{sc} V$)
of conics on $V$ is a smooth surface
(but possibly reducible, e.g., for $g=6$), 
and called the {\em Fano surface} of conics on $V$.
Every general point of $\Gamma$ corresponds to a good conic on $V$
(of type $(0,0)$).
(Here we use the assumption that $\car k=0$.)
For the proof, we refer to \cite{Ceresa-Verra} for $g=2$,
\cite{Collino-Murre-Welters} for $g=3$,
and \cite[Proposition~4.2.5]{AG5} for $g \ge 4$
(whose proof works also for $g=4$).
Consequently, we have
\begin{lem}
  \label{lem:good conic}
  Every prime Fano $3$-fold $V$ contains a good conic $E$,
  i.e., a conic of type $(0,0)$. 
  There exists an open dense subset of the Fano surface
  of $V$, which parametrises good conics $E \subset V$.
\end{lem}

\subsection{Hilbert-flag schemes}
\label{subsec:flag scheme}

In this subsection, we recall some basic results
on Hilbert-flag schemes.
See \cite{Kleppe87,Sernesi} for the proofs.
Given a projective scheme $Z$ over $k$,
we denote by $\HF Z$ the {\em Hilbert-flag schemes} 
(or the {\em nested Hilbert scheme}) of $Z$,
which parametrises
all pairs $(X,Y)$ of closed subschemes $X$ and $Y$ of $Z$
satisfying $X \subset Y$.
There are two natural projections
$pr_i: \HF Z \rightarrow \Hilb Z$ ($i=1,2$)
to the Hilbert scheme $\Hilb Z$ of $Z$,
sending $(X,Y)$ to $[X]$ for $i=1$, 
and to $[Y]$ for $i=2$.
We denote by $N_{(X,Y)/Z}$ the normal sheaf of $(X,Y)$ in $Z$
(see \cite[4.5.2]{Sernesi} for its definition),
which is a sheaf on $Z$ with support contained in $Y$.
By definition,  as $\mathcal O_Z$-modules,
$N_{(X,Y)/Z}$ is isomorphic to a subsheaf of
the direct sum $N_{X/Z} \oplus N_{Y/Z}$
of the normal sheaves of $X$ and $Y$.
Moreover, there exists a natural cartesian square
\begin{equation}
  \label{diag:cartesian}
  \vcenter{
    \xymatrix{
      N_{(X,Y)/Z} \ar[d]_{\pi_1} \ar[r]^{\pi_2} \ar@{}[dr]|\square & 
      N_{Y/Z} \ar[d]_{|_X} \\
      N_{X/Z} \ar[r]^{\pi_{X/Y}} & N_{Y/Z}\big{\vert}_X \\
  }}
\end{equation}
of homomorphisms of sheaves on $Z$, 
by which $N_{(X,Y)/Z}$ is characterized,
where $|_X$ is the restriction of sheaves, 
$\pi_{X/Y}: N_{X/Z} \rightarrow N_{Y/Z}\big{\vert}_X$ 
is the natural projection of normal sheaves,
and $\pi_i$ ($i=1,2$) are induced
by the projections to direct summands.

Suppose now that
the two closed embeddings 
$X \hookrightarrow Y$ and $Y \hookrightarrow Z$ are 
{\em regular embeddings},
whose definition can be found in ~\cite[\S D.1]{Sernesi}.
Then by ~\cite[Proposition~4.5.3]{Sernesi},
$H^0(Z,N_{(X,Y)/Z})$ and $H^1(Z,N_{(X,Y)/Z})$ respectively
represent the tangent space and the obstruction space
of $\HF Z$ at $(X,Y)$.
Moreover, 
it follows from a general theory that
\begin{equation}
  \label{ineq:dimension of flag}
  h^0(Z,N_{(X,Y)/Z})-h^1(Z,N_{(X,Y)/Z}) \le \dim_{(X,Y)} \HF Z
  \le h^0(Z,N_{(X,Y)/Z}),
\end{equation}
and $\HF Z$ is nonsingular at $(X,Y)$ if and only if
$\dim_{(X,Y)} \HF Z = h^0(Z,N_{(X,Y)/Z})$ (cf.~\cite[Lemma 7]{Kleppe87}).
The induced map $p_i=H^0(Z,\pi_i)$ by $\pi_i$
on the space of global sections
is the tangent map of $pr_i$ for each $i=1,2$.
Then the diagram \eqref{diag:cartesian} extends to
a commutative diagram 
\begin{equation}
  \label{diag:fundamental sequences of normal sheaves}
  \begin{CD}
    @. @. 0 @. 0 @.   \\
    @. @. @VVV @VVV  \\
    @. @. \mathcal I_{X/Y}\otimes_Y N_{Y/Z} @= \mathcal I_{X/Y}\otimes_Y N_{Y/Z} @.   \\
    @. @. @VVV @VVV  \\
    0 @>>> N_{X/Y} @>>> N_{(X,Y)/Z} @>{\pi_2}>> N_{Y/Z} @>>> 0 \\
    @. \Vert @. @V{\pi_1}VV @V{|_X}VV @. \\
    0 @>>> N_{X/Y} @>>> N_{X/Z} @>{\pi_{X/Y}}>> N_{Y/Z}\big{\vert}_X @>>> 0 \\
    @. @. @VVV @VVV  \\
    @. @. 0 @. 0   @. 
  \end{CD}
\end{equation}
of exact sequences of sheaves of $\mathcal O_Z$-modules.
Taking a long exact sequence,
we deduce from the first column
of \eqref{diag:fundamental sequences of normal sheaves}
the fundamental exact sequence 
\begin{equation}
  \label{seq:flag to hilb1}
  \begin{CD}
    0 @>>> H^0(Y,\mathcal I_{X/Y}\otimes_Y N_{Y/Z}) 
    @>>> H^0(Z,N_{(X,Y)/Z})
    @>{p_1}>> H^0(X,N_{X/Z})  \\
    @>>> H^1(Y,\mathcal I_{X/Y}\otimes_Y N_{Y/Z}) 
    @>>> H^1(Z,N_{(X,Y)/Z})
    @>{o_1}>> H^1(X,N_{X/Z})  \\
    @>>> H^2(Y,\mathcal I_{X/Y}\otimes_Y N_{Y/Z}) 
    @>>> \dots
  \end{CD}
\end{equation} 
of cohomology groups, and similarly from the first row
of \eqref{diag:fundamental sequences of normal sheaves}
the exact sequence
\begin{equation}
  \label{seq:flag to hilb2}
  \begin{CD}
    0 @>>> H^0(X,N_{X/Y}) 
    @>>> H^0(Z,N_{(X,Y)/Z})
    @>{p_2}>> H^0(Y,N_{Y/Z})  \\
    @>>> H^1(X,N_{X/Y}) 
    @>>> H^1(Z,N_{(X,Y)/Z})
    @>{o_2}>> H^1(Y,N_{Y/Z})  \\
    @>>> H^2(X,N_{X/Y}) 
    @>>> \dots
  \end{CD}
\end{equation}
(cf.~\cite[(2.7)]{Kleppe87}).
Here $o_i$ represents the maps on the obstruction spaces
induced by $pr_i$ for $i=1,2$.
By \cite[Lemma~A10]{Kleppe87}, if $H^1(Y,\mathcal I_{X/Y}\otimes_Y N_{Y/Z})=0$, 
then the morphism $pr_1$ is smooth at $(X,Y)$
(see also \cite[Theorem~1.3.4]{Kleppe81}).
Similarly one can deduce from \eqref{seq:flag to hilb2}
the fact that if $H^1(X,N_{X/Y})=0$, then $pr_2$ is smooth at $(X,Y)$
(cf.~\cite[Proposition 1.3.7]{Kleppe81}).
\begin{lem}[{cf.~\cite{Kleppe87}}]
  \begin{enumerate}
    \item If $H^1(Y,N_{Y/Z})=0$, then
    \[
      H^1(Z,N_{(X,Y)/Z})\simeq \coker \alpha_{X/Y},
    \]
    where $\alpha_{X/Y}$ is the composition
    $
    \alpha_{X/Y}: H^0(Y,N_{Y/Z}) 
    \overset{{\vert}_X}{\longrightarrow}
    H^0(X,N_{Y/Z}\big{\vert}_X)
    \overset{\partial_{X/Y}}{\longrightarrow}
    H^1(X,N_{X/Y})
    $
    of the restriction map $|_X$ and the coboundary map 
    $\partial_{X/Y}$ of the exact sequence in the second row of
    \eqref{diag:fundamental sequences of normal sheaves}.

    \item If $H^1(X,N_{X/Z})=0$, 
    then 
    \[
      H^1(Z,N_{(X,Y)/Z}) \simeq \coker \beta_{X/Y},
    \]
    where $\beta_{X/Y}$ is the composition
    $
    \beta_{X/Y}: H^0(X,N_{X/Z}) 
    \overset{\pi_{X/Y}}{\longrightarrow}
    H^0(Y,N_{Y/Z}\big{\vert}_X)
    \overset{\cup \mathbf k_X}{\longrightarrow}
    H^1(Y,\mathcal I_{X/Y}\otimes_Y N_{Y/Z})
    $
    of $\pi_{X/Y}$ and the coboundary map 
    $\cup \mathbf k_X$ of the exact sequence in the second column of
    \eqref{diag:fundamental sequences of normal sheaves}.
  \end{enumerate}
  \label{lem:obstruction space of flag}
\end{lem}
\Proof We show (1). 
By assumption, $H^1(Z,N_{(X,Y)/Z})$ is isomorphic to the cokernel
of the coboundary map 
$\partial_{X/Y}':H^0(Y,N_{Y/Z})\rightarrow H^1(X,N_{X/Y})$ 
of the exact sequence in the first row
of \eqref{diag:fundamental sequences of normal sheaves}.
By commutativity, we see that
$\partial_{X/Y}'$ factors through $\partial_{X/Y}$.
Similarly, (2) follows from the diagram.
\qed

\begin{rmk}
  \label{rmk:euler characteristic}
  The number in the left hand side of the inequality 
  \eqref{ineq:dimension of flag} is called
  the {\em expected dimension} of $\HF Z$ at $(X,Y)$.
  The Euler characteristic $\chi(Z,N_{(X,Y)/Z})$ of $N_{(X,Y)/Z}$ can be computed by
  the equation
  \[
    \chi(Z,N_{(X,Y)/Z})=\chi(X,N_{X/Y})+\chi(Y,N_{Y/Z})
    =\chi(Y,\mathcal I_{X/Y}\otimes_Y N_{Y/Z})+\chi(X,N_{X/Z}),
  \]
  which follows from the additivity on Euler characteristics and
  \eqref{diag:fundamental sequences of normal sheaves}.
  By definition, the support of $N_{(X,Y)/Z}$ is a closed subset of $Y$.
  Hence we have $H^i(Z,N_{(X,Y)/Z})=0$ for all integers $i > \dim Y$.
\end{rmk}

If $H^0(Y,\mathcal I_{X/Y}\otimes_{Y} N_{Y/Z})=0$,
then $pr_1$ induces an embedding of 
a neighborhood $U \subset \HF Z$ of $(X,Y)$ into $\Hilb Z$.
The following theorem is useful for estimating the (local) codimension
of the image of $\HF Z$ in $\Hilb Z$ at $[X]$,
and will be applied to Lemma~\ref{lem:codimension of S-maximal in hilb}.

\begin{thm}[cf.~\cite{Kleppe87,Kleppe89}]
  \label{thm:codimension of flag in hilb}
  Let $\mathcal O_{\Hilb Z,[X]}$ and 
  $\mathcal O_{\HF Z,(X,Y)}$ 
  denote the local rings
  of $\Hilb Z$ and $\HF Z$ at $[X]$ and $(X,Y)$, respectively.
  Let $\mathcal I_{X/Y}N_Y$ denote the sheaf 
  $\mathcal I_{X/Y}\otimes_{Y} N_{Y/Z}$ on $Y$, and 
  suppose that $H^0(Y,\mathcal I_{X/Y}N_Y)=0$.
  Suppose furthermore that $H^i(Z,N_{(X,Y)/Z})=0$ for $i=1,2$.
  Then we have
  \begin{equation}
    \label{ineq:codimension of HF in Hilb}
    h^1(Y,\mathcal I_{X/Y}N_Y)-h^2(Y,\mathcal I_{X/Y}N_Y)
    \le 
    \dim \mathcal O_{\Hilb Z,[X]} - \dim \mathcal O_{\HF Z,(X,Y)}
    \le
    h^1(Y,\mathcal I_{X/Y}N_Y).
  \end{equation}
  Here the inequality to the right is strict
  if and only if $\Hilb Z$ is singular at $[X]$.
  Let $\mathcal W_{X,Y}$ denote the unique irreducible
  component of $\HF Z$ passing through $(X,Y)$. Then,
  \begin{enumerate}
    \item If $h^1(Y,\mathcal I_{X/Y}N_Y)=0$ 
    or $h^2(Y,\mathcal I_{X/Y}N_Y)=0$, then
    $X$ is unobstructed in $Z$. Moreover,
    $\Hilb Z$ is generically smooth along $pr_1(\mathcal W_{X,Y})$.
    \item If $h^1(Y,\mathcal I_{X/Y}N_Y)=0$, 
    then $pr_1(\mathcal W_{X,Y})$ is an irreducible component of 
    $(\Hilb Z)_{\red}$.
    \item If $h^2(Y,\mathcal I_{X/Y}N_Y)=0$, 
    then $pr_1(\mathcal W_{X,Y})$ is of 
    codimension $h^1(Y,\mathcal I_{X/Y}N_Y)$ 
    in $\Hilb Z$ at $[X]$.
  \end{enumerate}
\end{thm}
\Proof 
We see that the proof of \cite[Theorem 10]{Kleppe87} 
works in our general setting, although
it is assumed there that $V=\mathbb P^3$, 
$Y$ is a smooth surfaces (of degree $s$)
and $X$ is a smooth curve (of degree $d >s^2$).
In fact, it follows from a general theory that
\begin{equation}
  \label{ineq:dimension of Hilb}
  h^0(X,N_{X/Z}) -h^1(X,N_{X/Z}) \le \dim \mathcal O_{\Hilb Z,[X]} 
  \le h^0(X,N_{X/Z}),
\end{equation}
and moreover, we have $\dim \mathcal O_{\Hilb Z,X} = h^0(X,N_{X/Z})$
if and only if $\Hilb Z$ is nonsingular at $[X]$.
Since $H^1(Z,N_{(X,Y)/Z})=0$, we have
$\dim \mathcal O_{\HF Z,(X,Y)}=h^0(Z,N_{(X,Y)/Z})$.
By subtracting this number from \eqref{ineq:dimension of Hilb}, 
we have
\begin{align*}
h^0(X,N_{X/Z})-h^0(Z,N_{(X,Y)/Z})-h^1(X,N_{X/Z})
& \le \dim \mathcal O_{\Hilb Z,X} -\dim \mathcal O_{\HF Z,(X,Y)} \\
& \le h^0(X,N_{X/Z})-h^0(Z,N_{(X,Y)/Z}). 
\end{align*}
It follows from the exact sequence \eqref{seq:flag to hilb1} and
$H^1(Z,N_{(X,Y)/Z})=H^0(Y,\mathcal I_{X/Y}N_Y)=0$ that
$h^0(X,N_{X/Z})-h^0(Z,N_{(X,Y)/Z})=h^1(Y,\mathcal I_{X/Y}N_Y)$.
Since we have $H^2(Z,N_{(X,Y)/Z})=0$ also, the same sequence shows that
$h^1(X,N_{X/Z})=h^2(Y,\mathcal I_{X/Y}N_Y)$.
Thus we have obtained \eqref{ineq:codimension of HF in Hilb}.
The first part of (1) is clear, because we have
$\dim \mathcal O_{\Hilb Z,X} -\dim \mathcal O_{\HF Z,(X,Y)}\ge 0$
by assumption,
while the last part of (1), (2) and (3)
follow from the upper semicontinuity
on cohomology groups.
\qed

\begin{rmk}
  \begin{enumerate}
    \item As is mentioned in \cite[Remark 11]{Kleppe87},
    by replacing the middle term
    \[
      \dim \mathcal O_{\Hilb Z,[X]} - \dim \mathcal O_{\HF Z,(X,Y)}
    \]
    in \eqref{ineq:codimension of HF in Hilb} with 
    \begin{equation}
      \label{num:codimension of the image of HF in Hilb}
      \dim \mathcal O_{\Hilb Z,[X]} - \dim \mathcal O_{\HF Z,(X,Y)}
      +h^0(Y,\mathcal I_{X/Y}N_Y),
    \end{equation}
    we can still prove Theorem~\ref{thm:codimension of flag in hilb} 
    without assuming that $H^0(Y,\mathcal I_{X/Y}N_Y)=0$.
    In fact, since $h^0(Y,\mathcal I_{X/Y}N_Y)$ 
    is greater than or equal to
    the dimension of the fiber at $[X]$ of the restriction
    $pr_1': \mathcal W_{X,Y} \rightarrow \Hilb Z$ of $pr_1$
    to $\mathcal W_{X,Y}$, the number
    \eqref{num:codimension of the image of HF in Hilb} is non-negative,
    and similarly we obtain all of (1), (2) and (3) in
    Theorem~\ref{thm:codimension of flag in hilb}.
    
    \item Theorem~\ref{thm:codimension of flag in hilb}~(2) 
    more directly follows 
    from the smoothness (or more precisely, the flatness) of 
    the morphism $pr_1$ at $(X,Y)$ 
    (cf.~\cite[Thm.~1.3.4 and Cor.~1.3.5]{Kleppe81}).
  \end{enumerate}
\end{rmk}

\subsubsection{Stably degenerate curves, $Y$-maximal family}
\label{subsubsec:stably and maximal}

From now on, we assume that
$\Hilb Z$ is nonsingular at $[Y]$, and $X$ is a smooth connected curve.
Let $W_Y$ denote the irreducible component of $\Hilb Z$ 
passing through $[Y]$, and
let $Z \times W_Y \supset \mathcal Y \overset{pr_2}\longrightarrow W_Y$
be the universal subscheme over $W_Y$.
We consider the Hilbert scheme $\Hilb^{sc} \mathcal Y$
of smooth connected curves in $\mathcal Y$, 
i.e., the relative Hilbert scheme of $\mathcal Y/W_Y$.
Then $\Hilb^{sc} \mathcal Y$ is isomorphic to an open subscheme of 
the Hilbert-flag scheme $\HF^{sc} Z$ ($:=pr_1^{-1}(\Hilb^{sc} Z)$),
where $pr_1:\HF Z \rightarrow \Hilb Z$ is the first projection.
Let $pr_1': \Hilb^{sc} \mathcal Y \rightarrow \Hilb^{sc} Z$ denote 
the restriction of $pr_1$ to $\Hilb^{sc} \mathcal Y \subset \HF Z$.
\begin{dfn}
  \label{dfn:stably degenerate}
  $X$ is {\em stably ($Y$-)degenerate} if $pr_1'$ 
  is surjective in a (Zariski) open neighborhood of $[X] \in \Hilb^{sc} Z$.
\end{dfn}
By definition, $X$ is stably degenerate if and only if
there exists an open neighborhood $U_X \subset \Hilb^{sc} Z$ of $[X]$
such that for any member $X'$ of $U_X$,
there exists a deformation of $Y'$ of $Y$ in $Z$ such that
$X' \subset Y'$ and $[Y'] \in W_Y$.
The following is one of the most
fundamental results on the stability of degenerate curves.
\begin{lem}
  \label{lem:general principle}
  If $\HF^{sc} Z$ is nonsingular at $(X,Y)$
  and $pr_1$ is smooth at $(X,Y)$,
  then $X$ is stably degenerate and unobstructed in $Z$.
\end{lem}
\Proof
By a property of smooth morphisms,
$\Hilb^{sc} Z$ is nonsingular at $[X]$.
Since the smoothness is a local property, 
$pr_1'$ is smooth at $(X,Y)$ and $\Hilb^{sc} \mathcal Y$ is
nonsingular at $(X,Y)$.
Let $\mathcal W_{X,Y}$ denote the unique irreducible component of
$\Hilb^{sc} \mathcal Y$ passing through $(X,Y)$.
Then by the smoothness,
its image $W_{X,Y}:=pr_1'(\mathcal W_{X,Y})$ is 
an irreducible component of $\Hilb^{sc} Z$, 
and this is the only one passing through $[X]$.
Hence $pr_1'$ is dominant near $[X]$.
\qed

\medskip

In general, the image $W_{X,Y}$ of an irreducible component
$\mathcal W_{X,Y} \subset \Hilb^{sc} \mathcal Y$ is
just an irreducible closed subset
of $\Hilb^{sc} Z$, and 
called the {\em $Y$-maximal family of curves} (containing $X$)
(cf.~\cite{Mukai-Nasu,Nasu5}).

\subsection{Curves and $K3$ surfaces in a Fano $3$-fold}
\label{subsec:curve and k3}

In this section, we recall some results from \cite{Nasu5},
concerned with the deformations of curves and $K3$ surfaces
in a smooth Fano $3$-fold.
Lemmas~\ref{lem:non-lifting of curves} and \ref{lem:k3 and fano}
are two key lemmas to prove Theorems~\ref{thm:main1} and \ref{thm:main2}.

Let $V$ be a smooth Fano $3$-fold,
$S$ a smooth member of $|-K_V|$, i.e., a smooth $K3$ surface,
$C$ a smooth connected curve on $S$.
Since $K_S \sim 0$, 
we have by adjunction that $N_{C/S}\simeq K_C$ and
$N_{S/V}\simeq -K_V\big{\vert}_S$,
and then $H^1(C,N_{C/S})\simeq k$.
By the ampleness of $-K_V$, we see that
$H^i(S,N_{S/V})=0$ for all integers $i>0$,
and hence $S$ is unobstructed in $V$.  
Then by Lemma~\ref{lem:nonsingularity of flag} below,
the Hilbert-flag scheme $\HF^{sc} V$ of $V$ is nonsingular 
at $(C,S)$ of expected dimension $\chi(V,N_{(C,S)/V})$,
i.e., we have $H^1(V,N_{(C,S)/V})=0$,
if and only if there exists a first order deformation
$\tilde S$ of $S$ in $V$ to which $C$ does not lift.
By Remark~\ref{rmk:euler characteristic}, or 
more directly from \cite[Lemma 2.10]{Nasu5},
the Euler characteristic of $N_{(C,S)/V}$ is computed as
\begin{equation}
  \label{eqn:euler characteristic}
  \chi(V,N_{(C,S)/V})=(-K_V)^3/2+g(C)+1,
\end{equation}
where $g(C)$ denotes the arithmetic genus of $C$.
\begin{lem}[{cf.~\cite{Nasu5}}]
  \label{lem:nonsingularity of flag}
  Let $V$ be a smooth projective scheme,
  $S \subset V$ a smooth surface,
  $C \subset S$ a smooth curve. 
  \begin{enumerate}
    \item Suppose that $H^1(S,N_{S/V})=0$ and $H^1(C,N_{C/S})\simeq k$.
    Then $H^1(V,N_{(C,S)/V})=0$ if and only if
    there exists a first order deformation $\tilde S$ of $S$ in $V$,
    to which $C$ does not lift.
    \item If $H^2(S,N_{S/V})=0$, then $H^i(V,N_{(C,S)/V})=0$ for $i > 1$.
  \end{enumerate}
\end{lem}
\Proof 
(1) By \eqref{seq:flag to hilb2} and assumption,
there exists an exact sequence
\[
  \begin{CD}
    H^0(V,N_{(C,S)/V}) 
    @>{p_2}>> H^0(S,N_{S/V})
    @>>> H^1(C,N_{C/S})
    @>>> H^1(V,N_{(C,S)/V}) 
    \longrightarrow 0.
  \end{CD}
\]
Then since $H^1(C,N_{C/S})$ is of dimension $1$,
we have $H^1(V,N_{(C,S)/V})=0$ if and only if
$p_2$ is not surjective. Thus we have proved (1).
By Remark~\ref{rmk:euler characteristic}, 
it suffices to show $H^2(V,N_{(C,S)/V})=0$ for (2),
which follows from \eqref{seq:flag to hilb2} and that $\dim C=1$.
\qed

\medskip

The following example shows that 
even if the second projection $pr_2: \HF^{sc} V \rightarrow \Hilb^{sc} V$
is not surjective in any neighborhood of $[S]$,
its tangent map
$p_2: H^0(V,N_{(C,S)/V}) \rightarrow H^0(S,N_{S/V})$ at $(C,S)$
can be surjective
(and hence $H^1(V,N_{(C,S)/V})\ne 0$).

\begin{ex}
  \label{ex:obstruction of flag including lines}
  Suppose that $V$ is a prime Fano $3$-fold of genus $g$,
  and let $E$ be a line on $V$
  (i.e.~$E \simeq \mathbb P^1$ and $(-K_V.E)_V=1$).
  Then $E$ is called a {\em good line} on $V$, 
  if $N_{E/V}$ is of type $(0,-1)$, 
  and a {\em bad line}, otherwise, i.e., $N_{E/V}$ is of type $(1,-2)$.
  Suppose that $S$ contains $E$.
  Then we have
  \[
    H^1(V,N_{(E,S)/V})\simeq 
    \begin{cases}
      0 & (\mbox{$E$: good}) \\
      k & (\mbox{$E$: bad})
    \end{cases}.
  \]
  Moreover, $\HF^{sc} V$ is nonsingular at $(E,S)$
  if and only if $E$ is good.
  In fact, the exact sequence 
  $0 \rightarrow N_{E/S} \rightarrow N_{E/V} 
  \rightarrow N_{S/V}\big{\vert}_E \rightarrow 0$
  on $E$ splits if and only if $E$ is bad.
  Suppose that $E$ is good.
  Then by $H^1(E,N_{E/V})=0$, the coboundary map
  $\partial_{E/S}: H^0(E,N_{S/V}\big{\vert}_E)
  \rightarrow H^1(E,N_{E/S}) (\simeq k)$ is surjective.
  Note that the restriction map
  $H^0(S,N_{S/V}) \overset{|_E}\longrightarrow H^0(E,N_{S/V}\big{\vert}_E)$
  is also surjective, because $E$ is a line.
  Then by Lemma~\ref{lem:obstruction space of flag},
  we have $H^1(V,N_{(E,S)/V})\simeq \coker \alpha_{E/S}=0$.
  Conversely, we suppose that $E$ is bad. Then by splitting,
  we have $\partial_{E/S}=0$ and 
  $H^1(V,N_{(E,S)/V})\simeq H^1(E,N_{E/S})\simeq k$.
  We also note that $pr_1$ is smooth at $(E,S)$ by $H^1(S,N_{S/V}(-E))=0$.
  Since $\Hilb^{sc} V$ is singular at $[E]$,
  so is $\HF^{sc} V$ at $(E,S)$.
  Since $E$ is not a complete intersection in $S$,
  i.e., $E \not\sim n(-K_V\big{\vert}_S)$ for any $n \in \mathbb Z$,
  we see that $pr_2$ is not surjective 
  in any neighborhood of $[S]$ by \cite{Moisezon},
  although $p_2$ is surjective for $(E,S)$.
\end{ex}

Repeating the same argument in
Example~\ref{ex:obstruction of flag including lines} for a conic $E$ on $V$, 
we obtain the following lemma.

\begin{lem}
  \label{lem:non-lifting of good lines and good conics}
  Suppose that $V$ is prime,
  and $E$ is a line or conic on $V$ contained in $S$.
  If $E$ is good, or a conic of type $(1,-1)$ (cf.~\eqref{eqn:type of conic}),
  then we have $H^1(V,N_{(E,S)/V})=0$.
\end{lem}

For the proof of Lemma~\ref{lem:non-lifting of curves} below,
we recall a criterion for the lifting of invertible sheaves.
Given a smooth projective scheme $X$ and 
an invertible sheaf $\mathcal L$ on $X$, 
we denote by $c(\mathcal L)$ 
the Atiyah extension class of $\mathcal L$ in
$H^1(X,\Omega_X)\simeq \Ext^1(T_X,\mathcal O_X)$.
Here $c(\mathcal L)$ is the image of the class of $\mathcal L$
in $\Pic X \simeq H^1(X,\mathcal O_X^\times)$,
under the map induced by the map
$d\log: \mathcal O_X^\times \rightarrow \Omega^1_X$
taking logarithmic derivatives (cf.~\cite[V, Ex.~1.8]{Hartshorne}).
Given an element $\tau \in H^1(X,T_X)$, 
i.e., an abstract first order deformation $\tilde X$ of $X$,
$\mathcal L$ lifts to an invertible sheaf on 
$\tilde X$ if and only if
the cup product $\tau \cup c(\mathcal L)$
via the pairing
$H^1(X,T_X) \times H^1(X,\Omega_X) 
\overset{\cup}\longrightarrow H^2(X,\mathcal O_X)$
is zero (cf.~\cite[Ex.~10.6]{Hartshorne10}).

\begin{lem}
  \label{lem:lifting of invertible sheaves}
  Let $V$ be a smooth projective scheme,
  $S$ a smooth closed subscheme of $V$,
  $\mathcal L$ an invertible sheaf on $V$, and
  $\mathcal M:=\mathcal L\big{\vert}_S$ its restriction to $S$.
  If $\tau \in H^1(S,T_S)$ is contained in the image
  of the coboundary map 
  $\delta: H^0(S,N_{S/V})\rightarrow H^1(S,T_S)$
  of the exact sequence
  \begin{equation}
    \label{ses:tangent and normal bundles}
    \begin{CD}
      0@>>> T_S @>{\iota}>> T_V\big{\vert}_S
      @>>> N_{S/V} @>>> 0,
    \end{CD}
  \end{equation}
  then we have $\tau \cup c(\mathcal M)=0$ in $H^2(S,\mathcal O_S)$.
\end{lem}
\Proof
Let $\iota$ be the inclusion of sheaves in 
\eqref{ses:tangent and normal bundles}.
Then its dual $\iota^{\vee}$ induces a map
$H^1(\iota^{\vee}): H^1(S,\Omega_V\big{\vert}_S) \rightarrow H^1(S,\Omega_S)$
on the cohomology groups.
There exists a commutative diagram
\[
  \xymatrix{
    H^1(V,\mathcal O_V^\times) \ar[d]_{|_S} \ar[rr]^{c}  & &
    H^1(V,\Omega_V) \ar[d]_{|_S} \\
    H^1(S,\mathcal O_S^\times) \ar[r]^{c} & H^1(S,\Omega_S) & \ar[l]_{H^1(\iota^\vee)}H^1(S,\Omega_V\big{\vert}_S). \\
  }
\]
Therefore, we have 
$c(\mathcal M)=H^1(\iota^{\vee})(c(\mathcal L)\big{\vert}_S)$
in $H^1(S,\Omega_S)$.
Then
\[
  \tau \cup c(\mathcal M)
  =\tau \cup H^1(\iota^{\vee})(c(\mathcal L)\big{\vert}_S)
  =H^1(\iota)(\tau) \cup c(\mathcal L)\big{\vert}_S,
\]
where the map 
$H^1(\iota): H^1(S,T_S)\rightarrow H^1(S,T_V\big{\vert}_S)$ 
is induced by $\iota$.
By assumption, we have $H^1(\iota)(\tau)=0$,
and hence we have finished the proof.
\qed

\begin{lem}[$\car k\ne 0$]
  \label{lem:non-lifting of curves}
  Let $i: S \hookrightarrow V$ be the closed embedding of $S$ into $V$,
  and $E$ an effective Cartier divisor on $S$ satisfying
  $
  H^1(S,\mathcal O_S(E))=H^1(V,N_{(E,S)/V})=0
  $.
  If the class of $C-bE$ in $\Pic S$ is contained in
  the image of the pullback map $i^*: \Pic V \rightarrow \Pic S$
  for some integer $b\ne 0$,
  then we have $H^1(V,N_{(C,S)/V})=0$.
\end{lem}
\Proof
By Lemma~\ref{lem:nonsingularity of flag},
there exists a first order deformation $\tilde S$ of $S$ in $V$
to which $E$ does not lift.
Let $\alpha \in H^0(S,N_{S/V})$ be the global section corresponding to $\tilde S$,
and $\tau:=\delta(\alpha)$ its image in $H^1(S,T_S)$
by the coboundary map $\delta$ of \eqref{ses:tangent and normal bundles}.
Then by Lemma~\ref{lem:lifting of invertible sheaves},
we have $\tau \cup c(\mathcal O_S(C-bE))=0$ in $H^2(S,\mathcal O_S)$.
On the other hand, since $c: H^1(S,\mathcal O_S^{\times}) \rightarrow
H^1(S,\Omega_S)$ is a group homomorphism, we have
$\tau \cap c(\mathcal O_S(C))=b \tau \cap c(\mathcal O_S(E))$.
If the invertible sheaf $\mathcal O_S(E)$ on $S$ lifts to $\tilde S$, 
then so does $E$
as a closed subscheme of $S$ by $H^1(S,\mathcal O_S(E))=0$
(cf.~\cite[Remark~4.5]{TDTEGA5}).
Therefore, we have 
$\tau \cap c(\mathcal O_S(E))\ne 0$ and hence
$\tau \cap c(\mathcal O_S(C))\ne 0$ by $b \ne 0$.
Then $\mathcal O_S(C)$ does not lift to $\tilde S$,
hence neither does $C$ (cf.~\cite[Ex.~6.7]{Hartshorne10}).
Thus we have finished the proof 
by Lemma~\ref{lem:nonsingularity of flag} again.
\qed

\medskip

We define a Cartier divisor $D$ on $S$ by
\begin{equation}
  \label{eqn:definition of D}
  D:=C+K_V\big{\vert}_S.
\end{equation}
Then since $\mathcal I_{C/S}\otimes_S N_{S/V}
\simeq \mathcal O_S(-D)$, the Serre duality shows
that
\begin{equation}
  \label{isom:cohomologies on D}
  H^i(S,\mathcal I_{C/S}\otimes_S N_{S/V})
  \simeq H^{2-i}(S,D)^{\vee}
\end{equation}
for all integers $i$.
Applying Theorem~\ref{thm:codimension of flag in hilb},
we have the following lemma.
\begin{lem}
  \label{lem:codimension of S-maximal in hilb}
  Let $W_{C,S} \subset \Hilb^{sc} V$
  be the $S$-maximal family of curves containing $C$.
  (See \S\ref{subsubsec:stably and maximal} for its definition.)
  Suppose that $H^0(S,-D)=H^1(V,N_{(C,S)/V})=0$.
  Then $\dim W_{C,S}=(-K_V)^3/2+g(C)+1$, and we have
  \[
    h^1(S,D)-h^0(S,D)
    \le \dim \mathcal O_{\Hilb^{sc} V, [C]} - \dim W_{C,S} 
    \le h^1(S,D).
  \]
  Moreover, $\Hilb^{sc} V$ is singular along $W_{C,S}$ if and only if
  the inequality to the right is strict.
  In particular,
  \begin{enumerate}
    \item If $h^1(S,D)=0$ or $h^0(S,D)=0$, then
    $C$ is unobstructed in $V$, and moreover
    $\Hilb^{sc} V$ is generically smooth along $W_{C,S}$.
    \item If $h^1(S,D)=0$, then $W_{C,S}$ is an irreducible component of 
    $(\Hilb^{sc} V)_{\red}$.
    In particular, $C$ is stably degenerate.
    \item If $h^0(S,D)=0$, then $W_{C,S}$ is of codimension $h^1(S,D)$
    in $\Hilb^{sc} V$.
  \end{enumerate}
\end{lem}
\Proof
Since $H^0(S,\mathcal I_{C/S}\otimes_S N_{S/V})\simeq H^0(S,-D)=0$, 
$\HF^{sc} V$ and $W_{C,S}$ are locally isomorphic
in a neighborhood of $(C,S)$,
and hence $\dim W_{C,S}=\dim \mathcal O_{\HF V,(C,S)}=(-K_V)^3/2+g(C)+1$.
\qed

\medskip

Finally we recall a result from \cite{Nasu5}. 
Given a curve $E$ on $S$, 
we define the {\em $\pi$-map} 
\begin{equation}
  \label{map:pi-map}
  \pi_{E/S}(E): 
  H^0(E,N_{E/V}(E))
  \longrightarrow
  H^0(E,N_{S/V}(E)\big{\vert}_E)
\end{equation}
for $E$ (or for the pair $(E,S)$) as the map on cohomology groups
induced by the sheaf homomorphism
$[N_{E/V} \overset{\pi_{E/S}}\longrightarrow N_{S/V}\big{\vert}_E] 
\otimes_E \mathcal O_E(E)$.
The following is another key lemma to prove 
Theorems~\ref{thm:main1} and \ref{thm:main2}.

\begin{lem}[{\cite[Theorem~1.2 and Corollary~1.3]{Nasu5}}]
  \label{lem:k3 and fano}
  Suppose that $H^1(V,N_{(C,S)/V})=0$.
  Let $D$ be the divisor on $S$ defined by \eqref{eqn:definition of D}, 
  and suppose that $D \ge 0$ and $D^2 \ge 0$.
  If there exists a $(-2)$-curve $E$ on $S$ such that
  \begin{enumerate}
    \renewcommand{\theenumi}{{\alph{enumi}}}
    \renewcommand{\labelenumi}{{\rm (\theenumi)}}
    \item $E.D=-2$,
    \label{item:intersection}
    \item $H^1(S,D-3E)=0$, and
    \label{item:vanishing}
    \item the $\pi$-map $\pi_{E/S}(E)$ is not surjective,
    \label{item:non-surjective}
  \end{enumerate}
  then $C$ is stably degenerate, and obstructed in $V$.
\end{lem}

\begin{rmk}
  \label{rmk:S-abnormal}
  If an effective divisor $D$ on a $K3$ surface $S$
  with $D^2\ge 0$ satisfy the above two conditions
  \eqref{item:intersection} and \eqref{item:vanishing}
  for some $(-2)$-curve $E$,
  then we have $h^1(S,D)=1$.
  (See \cite[Claim~4.1]{Nasu5} for the proof.)
  Thus in applying Lemma~\ref{lem:k3 and fano},
  the tangent map $p_1$ of the first projection $pr_1$
  is not surjective and $\coker p_1$ is of dimension $1$
  (cf.~\eqref{ses:flag to hilb1-2}).
\end{rmk}

\section{Non-reduced components of Hilbert schemes}
\label{sect:non-reduced}

In this section, we prove Theorem~\ref{thm:main1}.
As a result, we obtain an example
of a generically non-reduced component of
the Hilbert scheme $\Hilb^{sc} X$ (of Mumford type)
for every smooth Fano $3$-fold $X$ with $\Pic X \simeq \mathbb Z$ 
of any index $r$
(cf.~Corollary~\ref{cor:fano of picard number 1}, 
Table~\ref{table:GNRC of Mumford type}).

\subsection{Construction}
\label{subsec:construction}

Let $V$ be a prime Fano $3$-fold of genus $g$.
By Lemma~\ref{lem:good conic},
there exists a good conic $E\simeq \mathbb P^1$ on $V$.
Then we have $(-K_V.E)_V=2$ and $N_{E/V}$ is trivial.
As in the proof of \cite[Lemma~4.2.1]{AG5}, 
there exists a smooth member $S \in |-K_V|$ containing $E$.
Then $S$ is a $K3$ surface, and we have 
the self-intersection number $E^2=-2$ on $S$.
Put $\mathbf h:=-K_V\big{\vert}_S$, an ample divisor on $S$.
Then $\mathbf h^2=(-K_V\big{\vert}_S)^2=(-K_V)^3=2g-2$ and $\mathbf h.E=2$.
We consider a complete linear system
\[
  \Lambda:=|2\mathbf h+2E|
\]
of divisors on $S$.
By using the intersection numbers on $S$, we can show
that $\mathbf h+E$ is nef and big.
It is known that for every nef and big line bundle $L$ on a $K3$ surface,
$L^k$ is globally generated if $k \ge 2$
(cf.~\cite[Chap.~3, Remark 3.4]{Huybrechts16}).
Therefore, by Bertini's theorem,
$\Lambda$ contains a smooth connected curve $C$.
The degree and the genus of $C$ are computed as
$d(C)=C.\mathbf h=4g$ and $g(C)=C^2/2+1=4g+1$, respectively.

Let $W \subset \Hilb^{sc}_{4g,4g+1} V$ be the family of such curves 
$C \subset V$, i.e., 
smooth connected curves $C$ contained in a smooth $S \in |-K_V|$, 
and such that $C$ is a member of $\Lambda$ for some good conic $E$ on $V$.
For each $C \in W$, the surface $S$ and the conic $E$ 
are uniquely determined by $C$,
because we deduce
$h^0(V,\mathcal I_{C/V} \otimes_V \mathcal O_V(-K_V))\simeq k$
from the exact sequence
\begin{equation}
  \label{ses:ideal_CSV}
  \begin{CD}
    [0 @>>> \mathcal O_V(K_V)
      @>>> \mathcal I_{C/V} 
      @>{|_S}>> \mathcal O_S(-C)
      @>>> 0] \otimes_V \mathcal O_V(-K_V)
  \end{CD}
\end{equation}
and moreover, $E$ is recovered from $C$ and $S$ as
the unique base component of the linear system $|C+\mathbf h|$ on $S$.
Thus there exists a morphism
$W \rightarrow \HF^{sc} V$ to the Hilbert-flag scheme $\HF^{sc} V$ of $V$,
sending $[C]$ to $(E,S)$,
and the fiber at $(E,S)$ is isomorphic to an open subset of 
$\Lambda \simeq \mathbb P^{4g+1}=\mathbb P^{g(C)}
$
of smooth curves.
As in \S\ref{sect:introduction},
the pairs $(E,S)$ such that $E \subset S \subset V$
are parametrised by an open subset $U$
of $\mathbb P^{g-2}$-bundle over an open surface 
$\Gamma' \subset \Gamma$.
Thus there exists a diagram
\begin{equation}
  \label{diag:description of W}
  \vcenter{  
    \xymatrix{
      C \ar@{}[r]|*{\in} \ar@{|->}[d] & W^{(5g+1)} \ar@{}[r]|*{\subset} \ar[d]^{\mbox{\tiny $\mathbb P^{4g+1}$-bundle}}  & \Hilb^{sc}_{4g,4g+1} V \\
      (E,S) \ar@{}[r]|*{\in} \ar@{|->}[d] & U^{(g)} \ar@{}[r]|*{\subset} \ar[d]^{\mbox{\tiny $\mathbb P^{g-2}$-bundle}} & \HF^{sc} V \ar[d]^{} \\
      E \ar@{}[r]|*{\in} & \Gamma'^{(2)} \ar@{}[r]|*{\subset} & \Hilb^{sc}_{2,0} V
    }}
\end{equation}
of fiber bundles, 
where the upper script ${}^{(d)}$ of $X^{(d)}$ denotes $d=\dim X$.
Since $E$ is a good conic on $V$,
by Lemmas~\ref{lem:non-lifting of good lines and good conics} and
\ref{lem:non-lifting of curves},
we have $H^1(V,N_{(C,S)/V})=0$.
Then the Hilbert-flag scheme $\HF^{sc} V$ is nonsingular at $(C,S)$
of expected dimension
\[
  h^0(V,N_{(C,S)/V})=\chi(V,N_{(C,S)/V})=(g-1)+(4g+1)+1=5g+1
\]
by \eqref{eqn:euler characteristic}.
Let $D$ be the divisor on $S$ 
defined by \eqref{eqn:definition of D}.
Then it follows from
\eqref{seq:flag to hilb1} and \eqref{isom:cohomologies on D}
that there exists an exact sequence
\begin{equation}
  \label{ses:flag to hilb1-2}
  \begin{CD}
    0 @>>> H^0(V,N_{(C,S)/V}) @>{p_1}>> H^0(C,N_{C/V})
    @>>> H^1(S,D)^{\vee} @>>> 0.
  \end{CD}
\end{equation}
We will show that $H^1(S,D)\simeq k$, which implies 
that $h^0(C,N_{C/V})=h^0(V,N_{(C,S)/V})+1=5g+2$.
In fact, since $D\sim \mathbf h+2E$, 
we have $H^i(S,D-E)=0$ for $i=1,2$.
Then it follows from the exact sequence
\begin{equation}
  \label{ses:D and E}
  \begin{CD}
    0 @>>> \mathcal O_S(D-E) @>>> \mathcal O_S(D) @>>> \mathcal O_E(D) @>>> 0
  \end{CD}
\end{equation}
on $S$ that $H^1(S,D)\simeq H^1(E,\mathcal O_E(D))$.
Since $E.D=-2$ and $E\simeq \mathbb P^1$, we have $H^1(S,D)\simeq k$.
Thus we have a dichotomy between
\begin{enumerate}
  \renewcommand{\labelenumi}{{\rm (\Alph{enumi})}}
  \item The closure $\overline W$ of $W$ in $\Hilb^{sc} V$ 
  is an irreducible component of $(\Hilb^{sc} V)_{\red}$,
  and moreover $\Hilb^{sc} V$ is singular along $W$, and
  \item There exists an irreducible component $Z$ such that
  $\dim Z > \dim \overline W$ 
  and $\Hilb^{sc} V$ is generically nonsingular along $W$.
\end{enumerate}

\subsection{Proof of non-reducedness}
\label{subsec:non-reduced}

We prove that the case $(B)$ of the dichotomy in the previous subsection 
does not occur.

\paragraph{\bf Proof of Theorem~\ref{thm:main1}}
Let $W \subset \Hilb^{sc} V$ be the family of curves $C$ in $V$ 
as above. Then by the dichotomy,
it suffices to show that $C$ is obstructed in $V$.
We apply Lemma~\ref{lem:k3 and fano} to $C$.
It is easy to see that $D^2=2g-2>0$.
Moreover, since $N_{E/V}$ is globally generated,
the $\pi$-map $\pi_{E/S}(E)$ (cf.~\eqref{map:pi-map})
is not surjective by \cite[Lemma~2.14]{Nasu5}.
Thus we have only to check that $H^1(S,D-3E)=H^1(S,\mathbf h -E)=0$.
Since $E$ is a conic on $V \subset \mathbb P^{g+1}$, 
$E$ is linearly normal, i.e., 
$H^1(\mathbb P^{g+1},\mathcal I_{E/\mathbb P}(1))=0$.
Therefore the restriction map
$
H^0(V,-K_V) \rightarrow H^0(E,-K_V\big{\vert}_E)
$
to $E$ is surjective, and so is the map
$
H^0(S,\mathbf h) \rightarrow H^0(E,\mathbf h\big{\vert}_E).
$
Thus we deduce $H^1(S,\mathbf h-E)=0$ from 
\eqref{ses:D and E} with $D=\mathbf h$
and that $H^1(S,\mathbf h)=0$.
Thus we have finished the proof.
\qed

\medskip

\begin{rmk}
  By construction, the closure $\overline W$ of $W$ in $\Hilb^{sc} V$ is 
  nothing but the $S$-maximal family $W_{C,S}$
  of curves in $V$ containing $C$
  (cf.~\S\ref{subsubsec:stably and maximal}).
  In fact, we have $W \subset W_{C,S}$ and 
  $\dim W=W_{C,S}=5g+1$ ($=\dim \mathcal O_{\HF^{sc} V,(C,S)}$).
\end{rmk}

\begin{rmk}
By using the same construction and the same proof,
we can show that $\Hilb^{sc} V$ contains {\em infinitely many} generically 
non-reduced components (cf.~\cite[Example 5.8]{Nasu5}).
In fact, for every integer $n \ge 2$, we define 
a complete linear system $\Lambda_n$ on $S$ by
\[
  \Lambda_n:= |n\mathbf h+nE|.
\]
Then every general member $C_n$ is a smooth connected curve on $S$
of degree $2ng$ and genus $n^2g+1$.
Moreover, $C_n$ is parametrised by an irreducible locally closed subset
$W_n \subset \Hilb^{sc}_{2ng,n^2g+1} V$ of dimension $(n^2+1)g+1$
and we have $h^0(C_n,N_{C_n/V})=\dim W_n+1$.
The above argument (for $n=2$) works for $W_n$ in general.
Indeed, the generic member $C_n$ of $W_n$ is obstructed in $V$
(cf.~Proposition~\ref{prop:with conic}).
Thereby, the closure $\overline W_n$ of $W_n$ is an irreducible component of
$(\Hilb^{sc}_{2ng,n^2g+1} V)_{\red}$ and 
$\Hilb^{sc}_{2ng,n^2g+1} V$ is generically non-reduced along 
$\overline W_n$.
In particular, $\Hilb^{sc} V$ contains infinitely many generically 
non-reduced components.
\end{rmk}

\begin{rmk}
  \label{rmk:exotic fano}
  For some prime Fano $3$-folds $V$,
  the Hilbert scheme $\Hilb^{sc}_{1,0} V$ of lines on $V$
  is known to contain a generically non-reduced
  component (cf.~\cite[Proposition 4.2.2]{AG5}). 
  In the following two cases,
  $\Hilb^{sc}_{1,0} V$ contains a generically non-reduced component 
  $\Gamma$, whose general point corresponds to 
  a bad line on $V$, i.e., a line on $V$ of type $(1,-2)$
  (cf.~Example~\ref{ex:obstruction of flag including lines}):
  \begin{enumerate}
    \item $V$ is a smooth quartic hypersurface
    $V_4 \subset \mathbb P^4$ (i.e.~$g=3$),
    and there exists a ruled surface $R$ swept out by lines from $\Gamma$
    and $R$ is a cone over a smooth plane curve of degree $4$, e.g.,
    $V_4$ is the Fermat quartic 
    $V_{\rm Fer}=(\sum_{i=0}^4 x_i^2=0) \subset \mathbb P^4_{x_0,\dots,x_4}$.
    \item $R$ is swept out by projective tangent lines to some curve 
    $B  \subset V$, e.g., $V$ is a Mukai-Umemura $3$-fold
    $V_{22} \subset \mathbb P^{13}$ (cf.~\cite{Mukai-Umemura83}).
  \end{enumerate}
\end{rmk}

By the following examples, for every smooth Fano $3$-fold $X$ of index $4$ and $3$,
there exists a generically non-reduced component of $\Hilb^{sc} X$
(of Mumford type).
\begin{ex}[Mumford~\cite{Mumford}, $r=4$]
  \label{ex:mumford}
  Let $S \subset \mathbb P^3$ be a smooth cubic surface, 
  $E$ a $(-1)$-$\mathbb P^1$ on $S$ 
  and $C \subset S$ a smooth member
  of the linear system 
  $|-K_{\mathbb P^3}\big{\vert}_S+2E|\simeq \mathbb P^{37}$ on $S$.
  Then $C$ is of degree $14$ and genus $24$.
  Such $C$'s are parametrized by
  $W^{(56)} \subset \Hilb^{sc} \mathbb P^3$, which is an open subset
  of a $\mathbb P^{37}$-bundle over 
  $|\mathcal O_{\mathbb P^3}(3)|\simeq \mathbb P^{19}$.
  Then the closure $\overline W^{(56)}$
  is an irreducible component of $(\Hilb^{sc} \mathbb P^3)_{\red}$ and 
  $\Hilb^{sc} \mathbb P^3$ is everywhere non-reduced along $W^{56}$.
\end{ex}

\begin{ex}[$r=3$]
  \label{ex:non-reduced index 3}
  Let $Q$ be a smooth hyperquadric in $\mathbb P^4$, 
  and $S$ a smooth complete intersection of $Q$ with
  some other hyperquadric, i.e., $S \sim (-2/3)K_Q$.
  Let $\mathbf h \sim \mathcal O_S(1) \in \Pic S$ 
  be the class of hyperplane sections of $S$.
  Since $S$ is a del Pezzo surface of degree $4$
  (namely, $-K_S \simeq \mathbf h$ and $\mathbf h^2=4$), 
  $S$ is isomorphic to a blown up of $\mathbb P^2$ (at $5$ points).
  Thus there exists
  a $(-1)$-curve $E\simeq \mathbb P^1$ on $S$
  such that $E.\mathbf h=1$ (i.e.~a line $E$).
  We consider a complete linear system
  $|-K_Q\big{\vert}_S+2E|=|3\mathbf h+2E|$ on $S$.
  Then its general member $C$ is a smooth connected curve on $S$
  of degree $14$ and genus $16$.
  Since $N_{C/S}\simeq K_C(1)$ and $N_{S/Q}\simeq K_S(3)$,
  we have for all $i>0$ that $H^i(C,N_{C/S})=H^i(S,N_{S/Q})=0$
  and hence $H^i(Q,N_{(C,S)/Q})=0$ by \eqref{seq:flag to hilb2},
  which implies that
  $\HF Q$ is nonsingular at $(C,S)$ of expected dimension
  $\chi(C,N_{C/S})+\chi(S,N_{S/Q})=\chi(C,K_C(1))+\chi(S,-2K_S)=
  (d(C)+g(C)-1)+13=42$.
  Since $H^1(S,N_{S/Q}(-C))\simeq H^1(S,-\mathbf h-2E)\simeq k$,
  it follows from \eqref{seq:flag to hilb1} that $h^0(C,N_{C/Q})=43$.
  Then the $S$-maximal family $W_{C,S} \subset \Hilb^{sc} Q$ 
  of curves containing $C$ is a closed subset of $\Hilb^{sc} Q$
  of codimension $1$ in $H^0(C,N_{C/Q})$.
  Since for the generic member $C'$ of $W_{C,S}$,
  the line $E'$ on $S$ determined by $E'=\Bs|C'-2\mathbf h|$ 
  is a good line on $Q$ (i.e.~$N_{E'/Q}$ is of type $(1,0)$), 
  by using the same technique in \cite{Mukai-Nasu},
  we can show that $C'$ is obstructed in $Q$.
  Then by the same argument shown in
  \S\ref{subsec:construction} using a dichotomy,
  we see that $W_{C,S}$ is an irreducible component of 
  $(\Hilb^{sc} Q)_{\red}$,
  and hence $\Hilb^{sc} Q$ is generically non-reduced along $W_{C,S}$.
  Consequently, $\Hilb^{sc} Q$ contains a generically non-reduced
  component of Mumford type.
\end{ex}

For the index $2$ case, we refer to \cite{Mukai-Nasu,Nasu4}.
As a result, we have proved the following.

\begin{cor}
  \label{cor:fano of picard number 1}
  For every smooth Fano $3$-fold $X$ with $\Pic X\simeq \mathbb Z$,
  the Hilbert scheme $\Hilb^{sc} X$ of smooth connected curves on $X$
  contains a generically non-reduced component of Mumford type
  (cf.~Table\ref{table:GNRC of Mumford type}).
\end{cor}

\section{Deformations of degenerate curves}
\label{sect:stably degenerate}

In this section, we discuss the deformations of curves $C$ 
on a prime Fano $3$-fold $V$ and prove Theorem~\ref{thm:main2}.
We focus on the problem on determining
whether $C$ is stably degenerate or not
(cf.~Definition~\ref{dfn:stably degenerate})
for curves $C$ lying on a smooth anticanonical member $S$ of $V$.
In \cite{Nasu4}, the same problem was discussed 
for a smooth del Pezzo $3$-fold $V$, 
with curves lying on a smooth half-anticanonical member 
$S \in |-\frac 12K_V|$.

Let $V$ be a prime Fano $3$-fold of genus $g$,
and $C$ a smooth connected curve on $V$.
We suppose that there exists a smooth member $S \in |-K_V|$
containing $C$.
We recall that a curve $C$ on $S$
is called a {\em complete intersection in $S$}
if $C \sim n\mathbf h$ in $\Pic S$ for some integer $n$,
where $\mathbf h$ is the class
$-K_V\big{\vert}_S$
(of hyperplane sections of $S$).

\subsection{Complete intersection case}

We first consider the complete intersection case.
Given a positive integer $n$,
let $\mathcal W_n \subset \HF^{sc} V$ denote the (maximal) family
of pairs $(C,S)$ of a smooth member $S \in |-K_V|$
and a smooth curve $C \subset S$ 
satisfying $C \sim n \mathbf h$ in $\Pic S$.
Then the second projection $pr_2$ induces a morphism
$pr_2': \mathcal W_n \rightarrow 
|-K_V|\simeq \mathbb P^{g+1}$
to the Hilbert scheme of hyperplane sections of $V$.
It is clear that $pr_2'$ is dominant, and
its fiber at $[S]$ is isomorphic to 
an open subset of the projective space
$|\mathcal O_S(C)|\simeq |\mathcal O_S(n\mathbf h)|$.
Therefore, $\mathcal W_n$ is an irreducible locally closed subset
of $\HF^{sc} V$, and hence so is its image 
$W_n:=pr_1(\mathcal W_n)$ in $\Hilb^{sc} V$
by the first projection $pr_1$.
\begin{prop}
  \label{prop:complete intersection}
  Let $n \ge 1$ be an integer,
  and let $C$ and $W_n \subset \Hilb^{sc} V$ be as above.
  Then $C$ is unobstructed and stably degenerate in $V$.
  Moreover, the closure $\overline W_n$ of $W_n$
  is an irreducible component of $\Hilb^{sc} V$
  of dimension $(n^2+1)(g-1)+3$ for $n \ge 2$, and $2g$ for $n=1$.
\end{prop}
\Proof By using the Riemann-Roch theorem on $S$, we compute that
$\dim |\mathcal O_S(C)|=g(C)=(n\mathbf h)^2/2+1=n^2(g-1)+1$.
Hence
$\dim \mathcal W_n
=\dim |-K_V|+\dim |\mathcal O_S(C)|
=(n^2+1)(g-1)+3$.
Let $pr_1'$ denote the restriction of $pr_1$ to $\mathcal W_n$.
Then its fiber $pr_1'^{-1}([C])$ at $[C]$ 
is isomorphic to the linear system
$\Lambda:=|\mathcal I_{C/V}\otimes_V \mathcal O_V(-K_V)|$ on $V$.
It follows from \eqref{ses:ideal_CSV} and $H^1(V,\mathcal O_V)=0$ that
\[
  \dim \Lambda
  =h^0(S,-K_V\big{\vert}_S-C))=h^0(S,(1-n)\mathbf h),
\]
which is equal to $1$ if $n=1$, and $0$ otherwise.
Thus we obtain the dimension of $W_n$ as stated in the proposition.
It is also easy to see that $h^0(V,N_{C/V})=\dim W_n$.
In fact, there exists an exact sequence
$
0 \rightarrow K_C \rightarrow 
N_{C/V} \rightarrow \mathcal O_C(-K_V) \rightarrow 0
$
on $C$, because $N_{S/V}\simeq -K_V\big{\vert}_S$.
Since $C$ is a complete intersection in $S$, this sequence splits
and the restriction map
$H^0(S,-K_V\big{\vert}_S)\rightarrow H^0(C,-K_V\big{\vert}_C)$
to $C$ is surjective. Therefore, we compute that
$
h^0(C,N_{C/V})
=h^0(C,K_C)+h^0(C,-K_V\big{\vert}_C)
=\dim |\mathcal O_S(C)|+\dim |-K_V|-\dim \Lambda
=\dim W_n.
$
Since $\dim W_n \le \dim \mathcal O_{\Hilb^{sc} V,[C]}$,
we see by \eqref{ineq:dimension of Hilb} that
$C$ is unobstructed in $V$. 
Then it follows from $H^1(S,N_{S/V}(-C))=0$ that $pr_1'$ is smooth at $(C,S)$.
This implies that $\HF^{sc} V$ is nonsingular at $(C,S)$ and
hence $C$ is stably degenerate by Lemma~\ref{lem:general principle}.
\qed

\begin{rmk}
  We note that in Proposition~\ref{prop:complete intersection},
  $\HF^{sc} V$ is nonsingular, but 
  not of expected dimension at $(C,S)$,
  because $H^1(V,N_{(C,S)/V})\simeq k$.
  Moreover, $W_n$ coincides with 
  the $S$-maximal family $W_{C,S}$ containing $C$
  for every member $C$ of $W_n$
  and a member $S \in \Lambda$
  (cf.~\S\ref{subsubsec:stably and maximal}).
\end{rmk}

\subsection{Non-complete intersection case}
We next consider the case where $C$ is not a complete intersection in $S$.
For a technical reason (cf.~Remark~\ref{rmk:technical reason}),
we assume that the class of $C$ in $\Pic S$ is generated by $\mathbf h$
together with a line $E$ (cf.~Proposition~\ref{prop:with line}) 
or a conic $E$ (cf.~Proposition~\ref{prop:with conic})
with non-negative coefficients, i.e.,
$C \sim a\mathbf h +bE$ with two integers $a\ge 0$ and $b \ge 0$.
Since every line and every conic on $V$ 
is stably degenerate, and
there are many references, e.g., \cite{Iskovskih78,AG5},
for their (un)obstructedness in $V$,
in what follows, 
we also assume that $C$ is not a complete intersection in $S$
(i.e.~$b \ne 0$),
and $C \ne E$ (i.e.~$(a,b)\ne (0,1)$).
Let $W_{C,S} \subset \Hilb^{sc} V$ denote 
the $S$-maximal family of curves containing $C$.

\begin{prop}
  \label{prop:with line}
  Suppose that $E$ is a good line on $V$. Then
  \begin{enumerate}
    \item $C$ is stably degenerate, and unobstructed in $V$, and
    \label{item:stability and obstructedness}
    \item $W_{C,S}$ is an irreducible component of $(\Hilb^{sc} V)_{\red}$
    of dimension $(a^2+1)(g-1)+b(a-b)+2$.
    \label{item:S-maximal and dimension1}
  \end{enumerate}
\end{prop}

\Proof
We have $E.\mathbf h=1$ and $E^2=-2$.
Since $C\ne E$, we have $C.E\ge 0$, which implies that $a\ge 2b$.
By assumption, we have $b \ge 1$, and hence $H^1(V,N_{(C,S)/V})=0$ by 
Lemmas~\ref{lem:non-lifting of good lines and good conics}
and \ref{lem:non-lifting of curves}.
Let $D:=C+K_V\big{\vert}_S=(a-1)\mathbf h+bE$.
Since $D$ is effective and $D \not\sim 0$, we have $H^0(S,-D)=0$.
We note that $H^1(E,D\big{\vert}_E)=0$ by $D.E\ge -1$.
Since $D-E=(a-2b+1)\mathbf h+(b-1)(2\mathbf h+E)$ and 
$2\mathbf h+E$ is nef and big, we see that
$D-E$ is ample and hence $H^1(S,D-E)=0$.
Thus it follows from the exact sequence \eqref{ses:D and E}
that we have $H^1(S,D)=0$.
Then by Lemma~\ref{lem:codimension of S-maximal in hilb},
we obtain \eqref{item:stability and obstructedness}.
We compute the genus $g(C)$ of $C$ as
$g(C)=C^2/2+1=a^2(g-1)+ab-b^2+1$.
Thus by the same lemma,
we have proved \eqref{item:S-maximal and dimension1}.
\qed

\begin{prop}
  \label{prop:with conic}
  Suppose that $E$ is a good conic on $V$. Then
  \begin{enumerate}
    \item $C$ is stably degenerate,
    \label{item:stability}
    \item $C$ is obstructed in $V$ if and only if $a=b\ge 2$, and
    \label{item:obstructedness}
    \item $W_{C,S}$ is an irreducible component of $(\Hilb^{sc} V)_{\red}$
    of (the maximal) dimension $(a^2+1)g-(a-b)^2+1$
    (passing through $[C]$).
    \label{item:S-maximal and dimension2}
  \end{enumerate}
\end{prop}
\Proof 
The proof is similar to that of Proposition~\ref{prop:with line}.
We have $E.\mathbf h=2$ and $E^2=-2$.
Since $C\ne E$ and $b \ge 1$, we have $a \ge b$ and
$H^1(V,N_{(C,S)/V})=0$, respectively.
Put $D:=C+K_V\big{\vert}_S=(a-1)\mathbf h+bE$, a divisor on $S$,
as before.
Then $H^0(S,-D)=0$ by $D\ge 0$ and $D \ne 0$.
We consider the exact sequence \eqref{ses:D and E}
for the computation of $H^1(S,D)$.
First we note that $D-E=(a-b)\mathbf h+(b-1)(\mathbf h+E)$.
Since $\mathbf h+E$ is nef and big, so is $D-E$,
and hence we have $H^1(S,D-E)=0$
by the Kodaira-Ramanujam vanishing theorem.
Since $D.E=2(a-b)-2$, we see that
if $a > b$, then $H^1(E,D\big{\vert}_E)=0$
and hence $H^1(S,D)=0$ by \eqref{ses:D and E}.
Then in this case, 
by Lemma~\ref{lem:codimension of S-maximal in hilb},
we obtain the conclusions
\eqref{item:stability},
\eqref{item:obstructedness} and
\eqref{item:S-maximal and dimension2},
together with
$\dim W_{C,S}=(-K_V)^3/2+g(C)+1$, where $g(C)=a^2g-(a-b)^2+1$.
If $a=b=1$, then we have $H^1(S,D)=H^1(S,E)=0$, and hence the proof is done.

Suppose now that $a=b=n\ge 2$.
Then we have 
$D^2=(2g-2)(n-1)^2>0$,
$D\ge 0$,
$D.E=-2$ and
$D-3E=(n-1)\mathbf h+(n-3)E$.
Since we have $D-3E=\mathbf h-E$ for $n=2$,
and $D-3E=2\mathbf h+(n-3)(\mathbf h+E)$ for $n \ge 3$,
we conclude that $H^1(S,D-3E)=0$.
Since $E$ is a good conic on $V$,
the $\pi$-map $\pi_{E/S}(E)$ is not surjective,
as in the proof of Theorem~\ref{thm:main1}.
Thus by Lemma~\ref{lem:k3 and fano}, 
$C$ is stably degenerate and obstructed in $V$.
Moreover, we have $h^1(S,D)=1$ by Remark~\ref{rmk:S-abnormal},
and hence 
$h^0(C,N_{C/V})=h^0(V,N_{(C,S)/V})+1=
\dim W_{C,S}+1$ by \eqref{ses:flag to hilb1-2}.
Then we have the dichotomy between (A) and (B) 
in \S\ref{subsec:construction} for $W=W_{C,S}$ again.
Then the obstructedness immediately shows that
the case (B) does not occur.
Hence $W_{C,S}$ is an irreducible component of 
$(\Hilb_{sc} V)_{\red}$ of the maximal dimension
passing through $[C]$.
\qed

\medskip

Consequently, Theorem~\ref{thm:main2} has been proved
as a combination of Propositions~\ref{prop:complete intersection},
\ref{prop:with line} and \ref{prop:with conic}.

\begin{rmk}
  We have assumed that $(a,b)\ne (0,1)$, 
  i.e., $C$ is not $E$ itself,
  in Propositions \ref{prop:with line}, \ref{prop:with conic},
  and hence in Theorem~\ref{thm:main2}.
  If the condition
  [\ref{item:w.line}] or [\ref{item:w.conic}] of this theorem
  holds with $(a,b)=(0,1)$ (i.e.~$C=E$),
  then all the conclusions of the theorem are true, 
  except for the dimension of $\Hilb^{sc} V$ at $[C]$,
  which is computed as follows.
  Since $C=E$ is a good line or a good conic on $V$, we have $H^1(E,N_{E/V})=0$,
  and it follows from \eqref{ineq:dimension of Hilb} that
  $\dim_{[E]} \Hilb^{sc} V=\chi(E,N_{E/V})=(-K_V.E)_V$,
  that is equal to $1$ for [\ref{item:w.line}], 
  and $2$ for [\ref{item:w.conic}].
\end{rmk}

\begin{rmk}
  \label{rmk:technical reason}
  Here we restrict ourselves to the case where the class $[C]$ 
  of $C$ in $\Pic S$ is generated by $\mathbf h$ together with a line or a conic on $S$.
  Of course, this is not enough for determining 
  the stability of every $S$-degenerate non-complete intersection curves $C$ in $V$.
  It follows from a general principle 
  (cf.~Lemmas~\ref{lem:general principle} and \ref{lem:codimension of S-maximal in hilb})
  that if $H^1(V,N_{(C,S)/V})=H^1(S,D)=0$,
  then $C$ is stably degenerate and unobstructed in $V$.
  On the other hand, for some curves $C$ with $H^1(S,D)\ne 0$,
  it is difficult to prove that $C$ is stably degenerate.
  For example, if $[C]$ is generated by $\mathbf h$ and
  a $(-2)$-curve $E$ on $S$ of degree $E.\mathbf h>2$,
  then the intersection number $D.E=(C-\mathbf h).E$ can be less than $-2$
  (and then $h^1(S,D)>1$),
  and thereby Lemma~\ref{lem:k3 and fano} does not apply to $C$.
  For a detailed study of the stability of 
  $S$-degenerate curves on a smooth quartic $3$-fold $V$ ($g=3$),
  we refer to \cite{Nasu5},
  in which $S$ is assumed to be of Picard rank $2$.
\end{rmk}

\begin{rmk}
  Let $d(C)$ denotes the degree $(-K_V.C)_V$ of $C$.
  Then by counting dimensions, we see that
  if $H^1(V,N_{(C,S)/V})=0$ and $g(C)< d(C)-g$,
  then $C$ is not stably degenerate.
  In fact, by assumption we have
  \[
    \dim_{(C,S)} \HF^{sc} V=(-K_V)^3/2+g(C)+1=g(C)+g,
  \]
  while every component of the Hilbert scheme $\Hilb^{sc} V$
  is of dimension at least $\chi(C,N_{C/V})=d(C)$ 
  by \eqref{ineq:dimension of Hilb}.
  Thus if $pr_1': \Hilb^{sc} \mathcal S (\subset \HF^{sc} V)
  \rightarrow \Hilb V$ 
  is surjective in a neighborhood of $[C]$,
  then we have $g(C)+g\ge d(C)$.
\end{rmk}

\medskip

\paragraph{\bf Acknowledgments}
I should like to thank Prof.~Eiichi Sato 
for asking me a useful question, 
motivating me to research the subject of this paper.
I should like to thank Prof.~Shigeru Mukai
for useful comments on the tangent and obstruction spaces 
of Hilbert-flag schemes.
I should like to thank Prof.~Yuya Matsumoto for 
pointing out a mistake in an early version of the proof of
Lemma~\ref{lem:non-lifting of curves}.
I also thank the referee for the careful reading and valuable comments
improving the readability of this paper.
This work was supported in part by JSPS KAKENHI Grant Numbers JP17K05210, JP25220701.

\bibliography{mybib}
\bibliographystyle{abbrv}

\end{document}